\documentclass[12pt]{article}
\def\ignore#1{}
\usepackage{amsfonts}
\usepackage{amssymb}
\usepackage{amsmath}

\usepackage{epsfig, subfigure, verbatim}

\renewcommand{\emph}[1]{{\it #1}}
\newcounter{segcount}

\newenvironment{segment}
{\refstepcounter{segcount}\vspace{5mm}
\noindent{\bf \thesegcount. }}
{}

\newenvironment{statementnumbered}[3][]
{\refstepcounter{segcount}\vspace{5mm}
\noindent{\bf\thesegcount. #2}#1{\bf.}\ {\sl #3}}
{\nolinebreak[4] \nopagebreak[4] $\hfill \square$}

\newenvironment{statement}[3][]
{\vspace{5mm}\noindent{\bf #2}#1{\bf.} {\sl #3}}
{\nolinebreak[4] \nopagebreak[4] $\hfill \square$}

\newenvironment{statementnoboxnumbered}[3][]
{\refstepcounter{segcount}\vspace{5mm}
\noindent{\bf\thesegcount. #2}#1{\bf.} {\sl #3}}
{}

\newenvironment{statementnobox}[3][]
{\vspace{5mm}\noindent{\bf #2}#1{\bf.} {\sl #3}}
{}

\newenvironment{definitionnumbered}[2][]
{\refstepcounter{segcount}\vspace{5mm}
\noindent{\bf\thesegcount. #2}#1{\bf.}}
{}

\newenvironment{definition}[2][]
{\vspace{5mm}\noindent{\bf #2}#1{\bf.}}
{}

\newenvironment{resultnumbered}[3][]
{\refstepcounter{segcount}\vspace{5mm}
\noindent{\bf\thesegcount. #2}#1{\bf.} {\sl #3}
\vskip5mm\noindent {\bf Proof: }}
{\nopagebreak[4] $\hfill \square$}

\newenvironment{result}[3][]
{\vspace{5mm}
\noindent{\bf#2}#1{\bf.} {\sl #3}
{\\ \bf Proof: }}
{\nolinebreak[4] \nopagebreak[4] $\hfill \square$}

\newenvironment{risultnumbered}[3][]
{\refstepcounter{segcount}\vspace{5mm}
\noindent{\bf\thesegcount. #2}#1{\bf.} {\sl #3}
{\nopagebreak[4] \noindent \bf Proof: }}
{$\hfill \square$}

\newenvironment{risult}[3][]
{\vspace{5mm}
\noindent{\bf#2}#1{\bf.} {\sl #3}
{\noindent \bf Proof: }}
{\nolinebreak[4] \nopagebreak[4] $\hfill \square$}

\newenvironment{risultnoboxnumbered}[3][]
{\refstepcounter{segcount}\vspace{5mm}
\noindent{\bf\thesegcount. #2}#1{\bf.} {\sl #3}
{\nopagebreak[4] \noindent \bf Proof: }}
{}

\newenvironment{risultnobox}[3][]
{\vspace{5mm}
\noindent{\bf#2}#1{\bf.} {\sl #3}
{\nopagebreak[4] \noindent \bf Proof: }}
{}

\newenvironment{resultnoboxnumbered}[3][]
{\refstepcounter{segcount}\vspace{5mm}
\noindent{\bf\thesegcount. #2}#1{\bf.} {\sl #3}
{\\ \bf Proof: }}
{}

\newenvironment{resultnobox}[3][]
{\vspace{5mm}\noindent{\bf #2}#1{\bf.} {\sl #3}
{\\ \bf Proof: }}
{}

\newcommand{\seg}{\begin{segment}}
\newcommand{\segend}{\end{segment}}

\newcommand{\stmtnum}{\begin{statementnumbered}}
\newcommand{\stmtnumend}{\end{statementnumbered}}

\newcommand{\stmt}{\begin{statement}}
\newcommand{\stmtend}{\end{statement}}

\newcommand{\stmtnoboxnum}{\begin{statementnoboxnumbered}}
\newcommand{\stmtnoboxnumend}{\end{statementnoboxnumbered}}

\newcommand{\stmtnobox}{\begin{statementnobox}}
\newcommand{\stmtnoboxend}{\end{statementnobox}}

\newcommand{\defnnum}{\begin{definitionnumbered}}
\newcommand{\defnnumend}{\end{definitionnumbered}}

\newcommand{\defn}{\begin{definition}}
\newcommand{\defnend}{\end{definition}}

\newcommand{\resnum}{\begin{resultnumbered}}
\newcommand{\resnumend}{\end{resultnumbered}}

\newcommand{\res}{\begin{result}}
\newcommand{\resend}{\end{result}}

\newcommand{\risnum}{\begin{risultnumbered}}
\newcommand{\risnumend}{\end{risultnumbered}}

\newcommand{\ris}{\begin{risult}}
\newcommand{\risend}{\end{risult}}

\newcommand{\risnoboxnum}{\begin{risultnoboxnumbered}}
\newcommand{\risnoboxnumend}{\end{risultnoboxnumbered}}

\newcommand{\risnobox}{\begin{risultnobox}}
\newcommand{\risnoboxend}{\end{risultnobox}}

\newcommand{\resnoboxnum}{\begin{resultnoboxnumbered}}
\newcommand{\resnoboxnumend}{\end{resultnoboxnumbered}}

\newcommand{\resnobox}{\begin{resultnobox}}
\newcommand{\resnoboxend}{\end{resultnobox}}

\newcommand{\cF}{{\cal F}}
\newcommand{\cG}{{\cal G}}
\newcommand{\cI}{{\cal I}}

\newcommand{\cK}{{\cal K}}
\newcommand{\cL}{{\cal L}}
\newcommand{\cM}{{\cal M}}
\newcommand{\cO}{{\cal O}}
\newcommand{\cP}{{\cal P}}
\newcommand{\cQ}{{\cal Q}}

\newcommand{\Fs}{\cF_{\subseteq}}
\newcommand{\Fi}{\cF_{\leq}}

\newcommand{\la}{}
\newcommand{\ra}{{_{\leq}}}
\newcommand{\lb}{}
\newcommand{\rb}{{_{\subseteq}}}
\newcommand{\mL}{{\mathbb L}}
\newcommand{\mLa}{{\mathbb L}^a}
\newcommand{\mLi}{{\mathbb L}_{\leq}}
\newcommand{\mLai}{{\mathbb L}^a_{\leq}}
\newcommand{\mLc}{{\mathbb L}^c}

\newcommand{\mLic}{{\mathbb L}^c_{\leq}}
\newcommand{\mLidc}{{\mathbb L}^{dc}_{\leq}}

\def\cstar{${}$\circledast${}$}

\title{Factorisations and characterisations of induced-hereditary and compositive properties\footnote{The results presented here are part of the first author's Ph.D. thesis, written under the supervision of the second author. Jim Geelen suggested one of the main results of this paper.
}}

\author{
Alastair Farrugia%
\thanks{\noindent The first author's doctoral studies in Canada were fully funded by the
Canadian government through a Canadian Commonwealth Scholarship.}  
\and
R.\ Bruce Richter
\thanks{Research supported by NSERC.}
\\ {\it $\{$afarrugia, brichter$\}$@math.uwaterloo.ca}
\\ Dept. of Combinatorics \& Optimization
\\ University of Waterloo, Ontario, Canada, N2L 3G1
\and
Gabriel Semani\v{s}in
\thanks{Research supported in part by Slovak VEGA Grant 1/0424/03.}
\\ {\it semanisin@science.upjs.sk}
\\Institute of Mathematics
\\Faculty of Science, P.J. \v{S}af\'arik University
\\Jesenn\'a 5, 041 54 Ko\v{s}ice, Slovakia
}

\begin{document}
\maketitle
\begin{abstract}
A graph property (i.e., a set of graphs) is induced-hereditary or
additive if it is closed under taking induced-subgraphs or disjoint
unions. If $\cP$ and $\cQ$ are properties, the product $\cP \circ \cQ$
consists of all graphs $G$ for which there is a partition of the
vertex set of $G$ into (possibly empty) subsets $A$ and  $B$ 
with $G[A] \in \cP$ and $G[B] \in \cQ$. 
A property is reducible if it is the
product of two other properties, and irreducible otherwise. 

We completely describe the few reducible induced-hereditary properties that
have a unique factorisation into irreducibles.  Analogs
of compositive and additive induced-hereditary properties are
introduced and 
characterised in the style of Scheinerman [{\em Discrete
Math}.\ {\bf 55} (1985) 185--193].   One of these  provides an
alternative proof that an additive hereditary property 
factors into irreducible additive hereditary properties.  
\end{abstract} 

\section{Introduction}
Graph $k$-colouring can be viewed as a partitioning problem:  given a
graph $G$, can we partition its vertex set $V(G)$ into sets
$V_1,V_2,\dots,V_k$ so that each induced subgraph $G[V_i]$ is an
edgeless graph?  Allowing other possibilities than just edgeless
graphs introduces many new problems, such as: can we partition
$V(G)$ into two sets $V_1$ and $V_2$ so that $G[V_1]$ is a line graph
and $G[V_2]$ is perfect?  If $\cP_1,\cP_2,\dots,\cP_k$ are fixed
sets of graphs (also called {\em properties}), their {\em product} is the property
$\cP_1\circ\cP_2\circ\cdots\circ\cP_k$ consisting
of those graphs $G$ for which $V(G)$ partitions into $V_1,\dots,V_k$ so that, for each
$i$, $G[V_i]\in\cP_i$.  Each $V_i$ may be empty, so $\cP_i \subseteq \cP_1 \circ \cdots \circ \cP_k$.
A property is {\em reducible} if it is the
product of two other properties, and {\em irreducible} otherwise. 
A graph property is {\em induced-hereditary} or {\em additive} if it is closed under taking induced-subgraphs or disjoint unions. 
For a survey of induced-hereditary poperties we refer the reader to \cite{bobrfr}.

Scheinerman showed that intersection properties and induced-hereditary compositive properties are the same \cite{sch-1} (see~\cite{sch-4, sch-0, sch-2, sch-3} for related work). We give a Scheinerman-like characterisation of induced-hereditary disjoint compositive properties in Section \ref{sec-comp}. 

Although initially this characterisation was only intrinsically interesting, Jim Geelen
suggested that it could be used to prove the existence of a factorisation of an 
induced-hereditary disjoint compositive property into indecomposable 
induced-hereditary disjoint compositive properties.  As we discuss at the end of Section \ref{sec-alt-canonical}, this gives a simple proof of an important part \cite{uft-2} of another fundamental result \cite{uft-1,discussiones}, the unique factorisation theorem: if $\cP$ is an induced-hereditary disjoint compositive property, then there is a unique list of irreducible induced-hereditary disjoint compositive properties
$\cP_1,\dots,\cP_k$ such that $\cP=\cP_1\circ\cP_2\circ\cdots\circ\cP_k$. 

The introduction of induced-hereditary disjoint compositive properties is useful for
the study of the complexity of partitioning problems.  In particular, the general theory helps to show \cite{farr-comp,lozin} in many general cases that, if $\cP$ is additive induced-hereditary and $\cQ$ is coadditive induced-hereditary, then it is NP-hard to determine if a graph is in $\cP\circ\cQ$ (which is neither additive nor coadditive, but is disjoint compositive). 

In Section~\ref{sec-geq-her} we consider easy (non-)unique factorisation results for a natural class of properties. In Section \ref{sec-non-unique}, we completely describe all induced-hereditary properties that have a unique factorisation into arbitrary properties. Section~\ref{sec-infinite} is a short discussion of infinite graphs, while Section~\ref{sec-her-ind-her} explores briefly the relationship between forbidden subgraphs and forbidden induced-subgraphs.

Except in Section~\ref{sec-infinite}, we consider only finite, simple, unlabelled graphs.

\section{$\geq$-hereditary properties\label{sec-geq-her}}

To help the reader gain familiarity with some of the concepts, we provide in this section a simple factorisation theorem for a natural class of properties. We use $G\subseteq H$ and $G\le H$ to mean that $G$ is a subgraph of $H$ and $G$ is an induced-subgraph of $H$, respectively. A property $\cP$ is {\em $\geq$-hereditary\/} if it is closed under induced supergraphs; that is, if $G\geq H$ and $H\in \cP$, then $G\in\cP$. Induced-hereditary properties are defined similarly, using $\leq$ instead of $\geq$.
The class $\mL_{\geq}$ of all $\geq$-hereditary properties is closed under products, i.e., if $\cP,\cQ\in\mL_{\geq}$, then $\cP\circ\cQ\in\mL_{\geq}$.  Note that if $\cP\in\mL_{\geq}$, then $\cP$ is infinite and additive.

If $\cP$ is in $\mL_{\geq}$ and $\cP = \cP_1 \cup \cP_2$, then $\cP=\cP_1\circ\cP_2$; this is true also for coverings of $\cP$ by three or more sets. Thus,
for every integer $r > 1$, any property in $\mL_{\geq}$ 
has uncountably many factorisations into $r$ 
properties from $\mL_{\geq}$.
In particular, no $\ge$-he\-re\-di\-ta\-ry property is irreducible, even if we insist that the factors be in $\mL_{\geq}$. And yet we can specify a canonical factorisation quite easily. 

For any graph $G$, let $-G := \{H \mid G \nleq H\}$ and $+G := \{H \mid G \leq 
H\}$. Properties of the form $-G$ are {\em elementary\/} (cf. 
\cite{brown}) because every in\-du\-ced-he\-re\-di\-ta\-ry property 
$\cP$ can be expressed as $\cP = \bigcap\{-G \mid {G \in \Fi(P)} \}$, where $\Fi(P)$ denotes the set of graphs $G\notin \cP$ such that, for every $v\in V(G)$, $G-v\in\cP$. By a result of Berger~\cite{berger}, elementary properties are irreducible;  as we have just seen, this does not hold at all for $+G$, but properties of the form $+G$ are still special.

A $\ge$-he\-re\-di\-ta\-ry property $\cP$
is {\em primitive\/} if, in every factorisation
$\cP = \prod_{i \in I} \cP_i$ into $\ge$-he\-re\-di\-ta\-ry properties $\cP_i$,
there is an $i\in I$ 
such that $\cP = \cP_i$.  A factorisation $\cP = \prod_{i \in
I} \cP_i$ is {\em minimal\/} if 
there is no $I' \subset I$ such that $\cP = \prod_{i \in I'} \cP_i$.
The set $\min \cP$ consists of all (non-null) $\le$-minimal elements of $\cP$. 

\risnum{Proposition}
{Let $\cP$ be a $\geq$-hereditary property. Then $\cP$ is primitive iff, for some graph $G$, $\cP =+G$. Moreover, $\cP$ has a unique minimal factorisation into primitive $\geq$-hereditary properties:
\begin{equation}
\cP = \prod_{G \in \min{\cP}}\ +G . \label{prod-min}
\end{equation}
}
Since $\cP =\bigcup\ \{+G \mid G \in \min{\cP}\}$, (\ref{prod-min}) is a valid factorisation of $\cP$. Thus, if $\cP$ is primitive, then $|\min \cP| = 1$. Conversely, let $\min \cP = \{G\}$, and consider a factorisation $\cP_1 \circ \cP_2 \circ \cdots$ of $\cP$ into $\geq$-hereditary properties. For some $s \geq 0$, $G$ has a 
$(\cP_1, \ldots, \cP_s)$-partition $(V_1, \ldots, V_s)$, where some $V_i$'s may be empty. For all $i$, $G[V_i]$ is in $\cP_i \subseteq \cP$, so by minimality of $G$, either $V_i = \emptyset$ or $V_i = V(G)$.
Thus $G$ must be in, say, $\cP_1$. By $\geq$-heredity, $\cP_1 \supseteq +G = \cP$, so $\cP_1 = \cP$, and thus $\cP$ is primitive.

Now let $\cP$ be an arbitrary $\geq$-hereditary property, and let $\prod_{G \in S} +G$ be a factorisation into primitive $\geq$-hereditary properties. As argued above, each $G \in \min \cP$ must be contained in some property $+H$, where $H$ is in $S$; by minimality of $G$, we must have $G = H$. So $\min{\cP} \subseteq S$, and we must have equality for the factorisation to be minimal.
\risnumend


\section{Uniquely factorisable induced-hereditary \\ properties\label{sec-non-unique}}

The set $\mLai$ of additive induced-hereditary properties has a unique factorisation theorem \cite{uft-2, discussiones} (see~\cite{di_paper} for a generalisation to hypergraphs), 
but there are properties in the set $\mLi$ of induced-hereditary properties that are not uniquely factorisable over $\mLi$~\cite[Section 4]{uft-1}. 
Szigeti and Tuza~\cite[Prob. 4]{szi-t}\label{szi-t-1} asked whether an additive (induced-)hereditary property could have a factorisation where
the factors are not all additive and (induced-)hereditary. 
We show here that, in fact, practically no property in $\mLi$ is uniquely factorisable, when the factors are not required to be from $\mLi$. 


The in\-du\-ced-he\-re\-di\-ta\-ry property $\cO = \cO(\infty)$ consists of all finite edgeless graphs, while $\cK = \cK(\infty)$ consists of all finite complete graphs; for $s \in \mathbb{N}$, $\cO(s) := \{\overline{K}_r \mid 1\le r \leq s\}$ and 
${\cK}(s)  := \{K_r \mid 1\le r \leq s\}$.  For an arbitrary property $\cP$, $\la \cP \ra$ is the intersection of all induced-hereditary properties containing $\cP$.  If $G$ and $H$ are graphs, then $G\cup H$ is the vertex-disjoint union of $G$ and $H$, while $G+H$ is the {\em join\/} of $G$ and $H$, which is obtained from $G\cup H$ by adding all the edges joining a vertex in $G$ with a vertex in $H$.

\risnum{Theorem\label{uni=irr-O2}}
{Let $\cP$ be 
an in\-du\-ced-he\-re\-di\-ta\-ry property.  Then $\cP$  has a
unique factorisation if and only if one of the following occurs:
\begin{enumerate}
\item $\cP$ is irreducible;
\item $\cP = \cO(r) \circ \cO(s)$, where $r,s \leq \infty$;
\item $\cP = \cK(r) \circ \cK(s)$, where $r,s \leq \infty$;
\item $\cP = \cO(r) \circ \cK(s)$, where $r,s < \infty$.
\end{enumerate}
}
{\em Necessity\/}.  
We suppose $\cP$ is reducible, say $\cP=\cP_1\circ \cP_2$. 
If $\cP_1$ and $\cP_2$ are not both in\-du\-ced-he\-re\-di\-ta\-ry, then
it is easy to see that
$\la(\cP_1)\ra\circ \la(\cP_2)\ra$ is a different factorisation of $\cP$.
So we need only consider the case where $\cP_1$ and $\cP_2$ are both in $\mLi$.

Suppose first that $K_2$ and $\overline{K}_2$ are both in $\cP_1$.  Let
$\cP'_1=\cP_1\setminus\{K_1\}$.  We claim $\cP=\cP'_1\circ \cP_2$.
For if $G\in\cP$, let $(V_1,V_2)$ be a $(\cP_1,\cP_2)$-partition of
$G$.  If $|V_1|\ne 1$, then $(V_1,V_2)$ is a
$(\cP'_1,\cP_2)$-partition of $G$.  So suppose $|V_1|=1$.  If
$|V_2|=0$, then $(\emptyset,V_1)$ is a $(\cP'_1,\cP_2)$-partition of
$G$.  Otherwise, let $v\in V_2$.  Since both $K_2$ and $\overline{K}_2$
are in $\cP'_1$, and $\cP_2$ is in\-du\-ced-he\-re\-di\-ta\-ry,
$(V_1\cup\{v\},V_2\setminus\{v\})$ is a
$(\cP'_1,\cP_2)$-partition of $G$.  In all cases,
$G\in\cP'_1\circ\cP_2$.

Similarly, if $K_2$ and $\overline{K}_2$ are both in $\cP_2$, then $\cP$ does not
have a unique factorisation.  Now note that,
for any in\-du\-ced-he\-re\-di\-ta\-ry property $\cQ$, if $K_2\notin \cQ$, then
$\cQ\subseteq \cO$, while if $\overline{K}_2\notin \cQ$, then
$\cQ\subseteq {\cK}$.  

Moreover, since $\cP_1$ and $\cP_2$ are in\-du\-ced-he\-re\-di\-ta\-ry, we must be
in case (2), (3) or (4), except if $\cP_1 = \cO(r)$ and $\cP_2 = \cK$, or 
$\cP_1 = \cO$ and $\cP_2 = \cK(s)$.
Note that we may assume that $\cP_1\ne \{K_1\}$ and $\cP_2\ne\{K_1\}$,
as otherwise we have either (2) or (3).

Suppose $\cP_2={\cK}$.  We claim
$\cP=(\cP_1\setminus\{K_1\})\circ 
\cP_2$.  For let $G\in \cP$.  Then any partition of $V(G)$ into an
independent set and a clique works, unless the independent
set has size 1. If the independent set has size 1, and its vertex is
not joined to every vertex of the clique, then we can make the
independent set have size 2.  If the independent set has size 1 and is
joined to all the vertices in the clique, then $G$ is a clique, and so
is in $\cP_2$.  
Similarly, if $\cP_1=\cO$, then
$\cP=\cP_1\circ(\cP_2\setminus\{K_1\})$.  

{\em Sufficiency\/}.  The sufficiency of (1) is trivial.  
Let $\cP = \cO(r) \circ \cO(s)$ have some factorisation $\cQ_1 \circ \cQ_2$. Since $\cP$ contains only bipartite graphs, $\cQ_1$ and $\cQ_2$ contain only edgeless graphs.
Let $r', s'$ be positive integers, with $r' \leq r$ if $r$ is finite and $s' \leq s$ if $s$ is finite. $\overline{K}_{r'} + \overline{K}_{s'}$ is in $\cP$, and has a unique partition into two independent sets, so we must have $\overline{K}_{r'} \in \cQ_1$ and $\overline{K}_{s'} \in \cQ_2$ say. 

If $r \leq s$ are both finite, then $\overline{K}_{r'} + \overline{K}_{r'}$ shows that both factors contain $\overline{K}_{r'}$ for $r' \leq r$, while $\overline{K}_r +  \overline{K}_{s'}$ shows that one factor (and, clearly, only one) contains $\overline{K}_{s'}$ for all $r < s' \leq s$, so $\cO(r) \subseteq \cQ_1$, $\cO(s) \subseteq \cQ_2$, and we must clearly have equality.
If $r$ or $s$ is infinite, the proof is similar.

The sufficiency of (3) follows from that of (2) by complementation, so we consider $\cP = \cO(r) \circ \cK(s)$; if $r$ or $s$ is $1$, then we are in case (2) or (3), so we assume  $r,s \geq 2$. Let $\cQ_1 \circ \cQ_2$ be an arbitrary factorisation of $\cP$. If there is a
$K_a\in\cQ_1$ and a $K_b\in \cQ_2$ such that $a\ge 2$ and $b\ge 2$, then
$K_a\cup K_b\in \cP$, a contradiction.  Thus, we can assume $\cQ_1$
has no complete graph of size at least 2.

As $K_1\in \cO(r)$ and $K_2\in\cK(s)$, we have $K_3\in \cP$ and, therefore,
either $K_2$ or $K_3$ is in $\cQ_2$.  Let $b>1$ be such that
$K_b\in\cQ_2$.  Then, for every $G\in\cQ_1$, $G\cup K_b$ is in $\cQ_1\circ \cQ_2$. 
Let $(V_1,V_2)$ be an $(\cO(r),\cK(s))$-partition of $G\cup K_b$.  
Since $b>1$, at least one vertex from $K_b$
is in $V_2$.  It follows that no vertex from $G$ can be in $V_2$, so
$G\le (G \cup K_b)[V_1]\in\cO(r)$. Hence $\cQ_1\subseteq \cO(r)$. 

Suppose that $\cQ_1 = \{K_1\}$. Since ${\overline{K}_r} + K_s$ is in $\cP$, its $(\{K_1\}, \cQ_2)$ partition would imply that $\cQ_2$ has a graph containing $K_{s+1}$ or ${\overline{K}_r} + K_{s-1}$; thus $\cQ_1 \circ \cQ_2$ has a graph containing $K_{s+2}$, or a graph containing ${\overline{K}_r} + (K_{s-1} \cup K_1)$. Since we can check that neither of the last two graphs is in $\cO(r) \circ \cK(s)$, $\cQ_1$ must contain ${\overline{K}_a}$ for some $a \geq 2$.
Now for every $H \in \cQ_2$, ${\overline{K}_a} + H \in \cO(r) \circ \cK(s)$, and it follows that $\cQ_2 \subseteq \cK(s)$.

Let $r' \leq r$ be a positive integer.  Note that $G={\overline{K}_{r'}}+K_s$ is in $\cO(r)\circ\cK(s)$.  There are essentially only two
$(\cO,{\cK})$-partitions of $G$ and one of these uses
$K_{s+1}$, which is not in $\cQ_1$ or $\cQ_2$. 
Therefore, the only possible $(\cQ_1,\cQ_2)$-partition of
${\overline{K}_{r'}}+K_s$ shows ${\overline{K}_{r'}}\in\cQ_1$, and thus $\cQ_1=\cO(r)$.

Similarly, if $s' \leq s$ is a
positive integer, then ${\overline{K}_r}\cup
K_{s'}\in\cO(r)\circ\cK(s)$. There are essentially only two
$(\cO,{\cK})$-partitions of $G$, one of which uses $\overline
K_{r+1}$, which is not in $\cQ_1$.  Hence $K_{s'}\in\cQ_2$, so
$\cQ_2 = \cK(s)$, as required.
\risnumend

\section{Compositivity\label{sec-comp}}
In this section we give Scheinerman's characterisation of intersection properties as induced-hereditary compositive properties and also provide analogs for additive induced-hereditary and induced-hereditary disjoint compositive properties.
These are useful (as we demonstrate in Section \ref{sec-alt-canonical} for hereditary compositive properties) for finding a factorisation of such a property into other properties of the same type.

An {\em intersection graph\/} $I(F)$ is obtained from a finite collection $F$ of sets by making the vertices to be the sets in $F$, with adjacency determined by nonempty intersection.  If $\cF$ is a set of sets, then the {\em intersection property\/} of $\cF$ is $\cI(\cF) := \{I(F) \mid F \subseteq \cF, |F| \textrm{ finite}\}$.

A property $\cP$ is {\em induced-hereditary compositive\/} if it is induced-hereditary and, for every $G$ and $H$ in $\cP$, there is another graph $K$ in $\cP$ containing both $G$ and $H$ as induced subgraphs.  If $K$ may always be chosen so that $G$ and $H$ are disjoint induced subgraphs of $K$, then $\cP$ is {\em induced-hereditary disjoint compositive} (for short, {\em indiscompositive\/}).

A {\em generating set\/}\index{generating set} for $\cP$ is a set
$\cG$ such that $\la \cG \ra = \cP$, i.e., $\cG \subseteq \cP$ and
every graph in $\cP$ is an  
in\-du\-ced-sub\-graph of some $G \in \cG$. A property has a
generating set iff it is in\-du\-ced-he\-re\-di\-ta\-ry. A (possibly
finite) generating set is {\em ordered\/} if its elements can be
listed as $G_1 \leq G_2 \leq \cdots.$ 

A {\em generating graph\/}\index{generating graph} for $\cP$ is a
(finite or infinite) graph $H$ such that $\cP = \{G \mid G \leq H,
|V(G)| \textrm{ finite}\}$. 
A property need not have a unique generating graph, as shown by Bonato
and Tardif~\cite{bonato2}\label{bonato2-1} for $\cP := \{G \mid $
every component of $G$ is a path$\}$.
Let $S$ be an infinite set of
positive integers, and let $P_S$ be the disjoint union of paths of
length $s$, $s \in S$. There are uncountably many such generating graphs, and
they are all in\-du\-ced-sub\-graphs of each other.  
Note that $\cP$ also has generating graphs which are not contained
in each other (e.g. the two-way infinite path, and the disjoint union
of $3$ one-way infinite paths). 

For any $L\in\cP$, $\cG[L]$ is the set $\{G \in \cG \mid L \leq G \}$. 
\label{def-ind-her-contain-H-1}
Mih\'ok et al.\ introduced the set $\cG[L]$, and observed that A $\Rightarrow$ B in Theorem \ref{sch}; the rest of the theorem is due to Scheinerman. 
The equivalence of A and C is essentially proved also in~\cite[Thm.\ 2.3]{bebrmo}.
We reproduce its proof here, as it leads to characterisations of indiscompositive properties (Theorem~\ref{indiscomp-character}) and a simplification of the factorisation for $\mLa$ (the set of additive hereditary properties) (Theorem~\ref{factorn-new}).

\risnum[~\cite{uft-1,sch-1, sch-2}]{Theorem\label{sch}}
{For any property $\cP$, the following are equivalent:
\begin{itemize}
\item[A.] $\cP$ is in\-du\-ced-he\-re\-di\-ta\-ry compositive;
\item[B.] $\cP$ has a generating set; moreover, for any graph $L \in \cP$, and any generating set $\cG$ of $\cP$, $\cG[L]$ is also a generating set;
\item[C.] $\cP$ has a (finite or infinite) ordered generating set $H_1 \leq  H_2 \leq \cdots;$
\item[D.] $\cP$ has a generating graph $H$; and
\item[E.] $\cP$ is an intersection property.
\end{itemize}
}
(A $\Rightarrow$ B). Since $\cP$ is in\-du\-ced-he\-re\-di\-ta\-ry, it is itself a generating set for $\cP$. Now for any graph $G \in \cP$, compositivity implies that there is some graph $L'_G$ that contains both $L$ and $G$; since $\cG$ is a generating set there is some graph $L_G \in \cG$ that contains $L'_G$, and thus $G \leq L_G \in \cG[L]$.

(B $\Rightarrow$ C). Let $\cG_1 = \{G_1, G_2, G_3, \ldots\}$ be a generating set for $\cP$. Define $H_1 := G_1$; for $i \geq 1$,  $\cG_{i+1} := \cG_i[H_i]$ is also a generating set, so it contains some graph $H_{i+1}$ containing $H_i$ and $G_{i+1}$. Thus $H_1 \leq H_2 \leq H_3 \leq \cdots $ is an ordered generating set.

(C $\Rightarrow$ D). There is no loss of generality in assuming $|V(H_j)| = j$ for all $j$. If the ordered generating set is finite, say $H_1 \leq \cdots \leq H_r$, then take $H$ to be $H_r$. Otherwise, for each $j$ we can label $V(H_j) = \{v_1, v_2, \ldots, v_{j}\}$ so that the vertices $\{v_1, \ldots, v_{i}\}$ 
induce  $H_i$ (as a labeled graph) whenever $i \leq j$. So if $v_i$
and $v_j$ are adjacent in some $H_k$, then they are adjacent in all
subsequent $H_r$. We now let $H$ be the infinite graph  with vertices
$\{v_1, v_2, \ldots\}$, where $v_i$ and $v_j$ are adjacent iff they
are adjacent in some $H_k$. 
If $G \leq H$, with $v_j$ being the vertex of $V(G)$ with largest
index, then $G \leq H_j$ so that $G \in \cP$; conversely, if $G$ is in
$\cP$, then it is contained in some $H_j$, and therefore is contained
in $H[v_1, \ldots, v_j]$. 

(D $\Rightarrow$ E). Let $u \sim v$ denote that $u$ and $v$ are
adjacent. Given the graph $H$ with vertices $v_1, v_2, \ldots,$ define
the family of sets $\cF := \{S_i, i \geq 1\}$, where  
$S_i := \{\{v_i, v_j\} \mid v_i \sim v_j\}$. Then 
$S_i \cap S_j \not= \emptyset$ iff $v_i \sim v_j$.

(E $\Rightarrow$ A). If $G_1,G_2\in \cP$ and $\cF$ is a family of sets
whose corresponding intersection property is $\cP$, then, for $i=1,2$,
there is a finite subset $\cF_i$ of $\cF$, such that $G_i$ is the
intersection graph of $\cF_i$.  If $K\le G_1$, then $V(K)$ is some
subset $\cF_K$ of $\cF_1$ and $K$ is the intersection graph of
$\cF_K$.  Hence $\cP$ is induced hereditary.  Evidently, the
intersection graph of $\cF_1\cup \cF_2$ is in $\cP$ and contains both
$G_1$ and $G_2$ as induced subgraphs, so $\cP$ is compositive. 
\risnumend
\newline

%

We can generalise this theorem to indiscompositive properties, using
an infinite Ramsey Theorem. 

\stmt[{~\cite[Thm. 1.5]{graham}\label{graham-2}}]{Ramsey's
Theorem\label{inf-ramsey}} 
{If the edges of an infinite clique are coloured with finitely many
colours, then there is an infinite monochromatic clique.  
}
\stmtend
\newline

For a graph $L$, $2 \cstar L$ is the set of all graphs that consist of
two disjoint copies of $L$
with arbitrary edges between the two copies. 
For a set $\cG$, $\cG[2\cstar L]$ is the set of graphs in
$\cG$ which contain, as an induced-subgraph, some graph in $2 \cstar L$. 

Let $G$ be a graph with vertex-set $\{v_1, \ldots, v_k\}$. A {\em
uniform graph} $X$ for $G$ consists of infinitely many copies of $G$,
with a ``uniform'' set of edges between each pair. More precisely,
label the copies of $G$ as $G^1, G^2, \ldots;$ for all distinct $i,j
\geq 1$ label $V(G^i)$ as $v_1^i, \ldots, v_k^i$, and let $G^i:G^j$
denote the graph $X[v_1^i, \ldots, v_k^i, v_1^j, \ldots, v_k^j]$. Then
$X$ is {\em uniform for $G$\/} if: 

\begin{itemize}
\item[(*)] for each $i$ and $x$, $\varphi: v_x^i \mapsto v_x$ is an
isomorphism from $G^i$ to $G$ 
\item[(**)] for each $p < q$ and $x$, 
$\psi: \left\{\begin{tabular}{lll}
$v_x^p$ & $\mapsto$ & $v_x^1$ \\
$v_x^q$ & $\mapsto$ & $v_x^2$ 
\end{tabular}\right.$
is an isomorphism from $G^p:G^q$ to $G^1:G^2$
\end{itemize}

\risnum{Theorem\label{indiscomp-character}}
{For a property $\cP$ the following are equivalent:
\begin{itemize}
\item[A.] $\cP$ is indiscompositive;
\item[B.] $\cP$ has a generating set; moreover, for any graph $L \in
\cP$, and any generating set $\cG$ of $\cP$, $\cG[2\cstar L]$ is also
a generating set; 
\item[C.] $\cP$ has an ordered generating set $G_1 \leq G_2 \leq
\cdots$ such that every $G_i$ contains two disjoint copies of
$G_{i-1}$; and 
\item[D.] $\cP$ has a generating graph $H$ that contains two disjoint
copies of itself, i.e., $H \in 2 \cstar H$; moreover, for any
countable generating graph $G$, there is a generating graph $H \geq G$
such that $H \in 2 \cstar H$. 
\end{itemize}
}
A $\Rightarrow$ B and B $\Rightarrow$ C are proved as in
Theorem~\ref{sch}.  

(C $\Rightarrow$ A).  Since $\cP$ has a generating set, by Theorem
\ref{sch} it is induced-hereditary.  For $H_1, H_2\in \cP$, let
$j_1,j_2$ be such that $H_i\leq G_{j_i}$, $i=1,2$.  For
$j=\max\{j_1,j_2\}$, $H_1, H_2\leq G_j$.  Since $G_{j+1}$ contains
disjoint copies of $G_j$, it follows that $G_{j+1}$ contains disjoint
copies of $H_1$ and $H_2$. 

(A $\Rightarrow$ D). $\cP$ has a generating graph $G$ because it is
induced-hereditary compositive; by the construction in
Theorem~\ref{sch}.D, there is such a graph with countable vertex-set,
say $\{v_1, v_2, \ldots\}$. For each $k \geq 1$, let $G_k$ be $G[v_1,
v_2, \ldots, v_k]$.  

For each $k \geq 1$, we claim that there is a uniform graph for $G_k$. 
Let $G_{k,0}$ be the graph with no vertices. For each $\ell \geq 1$, we can find a graph $G_{k,\ell} \in \cP$ on $k\ell$ vertices containing $G_{k,\ell-1}$ and an $\ell^{\textrm{\scriptsize th}}$ disjoint copy of $G_k$. We label 
$V(G_{k,\ell})$ so that $v_1^1, \ldots, v_k^1, \ldots, v_1^{\ell-1}, \ldots, v_{k}^{\ell-1}$ give us $G_{k,\ell-1}$ (as a labeled graph), while $v_1^\ell, \ldots, v_{k}^\ell$ give us the same labeled graph as $G_k$ (when we ignore the superscripts).

Now since $k$ is finite, there are only finitely many (say $s_k$) configurations of edges that can be placed between two labeled copies of $G_k$. We colour the edges of an infinite clique $K$ with vertex set $\{t_1,t_2,\dots\}$ with (at most) $s_k$ colours so that, for $p < q$, the colour of arc  $(t_p, t_q)$ 
corresponds to the configuration between the $p^{\textrm{\scriptsize th}}$ and $q^{\textrm{\scriptsize th}}$ copies of $G_k$ in $G_{k,q}$. By Ramsey's Theorem, $K$ contains an infinite monochromatic clique, which corresponds to a uniform graph for $G_k$.
\newline

For each $k$, let $H_k^1, \ldots, H_k^{r_k}$ be the set of uniform graphs for $G_k$; note that $r_k \leq s_k$. Note that each such graph contains a uniform graph for $G_{k-1}$. 

Form a graph with vertex set $V_1 \cup V_2 \cup \cdots,$ where each $V_k$ contains $r_k$ vertices. The $p^{\textrm{\scriptsize th}}$ vertex of $V_k$ is joined to the $q^{\textrm{\scriptsize th}}$ vertex of $V_{k-1}$ if $H_k^p \geq H_{k-1}^q$; clearly, every vertex in $V_k$ has at least one neighbour in 
$V_{k-1}$. By K\"onig's Infinity Lemma \cite{diestel} there is an infinite path with exactly one vertex in each $V_i$, which corresponds to a sequence $H_1 \leq H_2 \leq \cdots$ where, for each $k$, $H_k$ is uniform for $G_k$. The nested union $H$ of the $H_i$ contains two copies of itself. Any finite subgraph $X$ of $H$ is contained in some $H_k$, and is thus in $\cP$; conversely, any graph in $\cP$ is contained in $G$ and, thus, in $H$. So $H$ is a generating graph for $\cP$ as required.

We omit the straightforward proof of (D $\Rightarrow$ A).
\risnumend
\newline

Additive in\-du\-ced-he\-re\-di\-ta\-ry properties also have similar characterisations. We use $kG$ to denote the disjoint union of $k$ copies of $G$. 
The proof is simpler than for Theorem \ref{indiscomp-character} and so is omitted.

\statementnoboxnumbered{Theorem\label{additive-character}}
{For any property $\cP$ the following are equivalent:
\begin{itemize}
\item[A.] $\cP$ is additive in\-du\-ced-he\-re\-di\-ta\-ry;
\item[B.] $\cP$ has a generating set; moreover, for any graph $L \in \cP$, and any generating set $\cG$ of $\cP$, $\cG[2L]$ is also a generating set;
\item[C.] $\cP$ has an infinite ordered generating set $H_1 \leq H_2 \leq \cdots$ such that $2H_i \leq H_{i+1}$ for all $i$; and
\item[D.] $\cP$ has a generating graph $H$ such that $2H \leq H$; moreover, for any generating graph $G$, there is a generating graph $H \geq G$ such that $2H \leq H$.\hfill$\square$
\end{itemize}
}
A property $\cP$ is {\em coadditive\/} if, for every $G,H\in\cP$, $G+H\in\cP$.
Additive and co-additive induced-hereditary properties are clearly indiscompositive. The set of split graphs is neither additive nor co-additive, though it is
the product of an additive and a  co-additive property. There are also irreducible properties that are indiscompositve but not additive or coadditive.  For example, let $K'$ be obtained from an
infinite clique $K$, by adding, for each $v\in V(K)$, a new vertex $v'$ adjacent only to
$v$. By Theorem~\ref{indiscomp-character}.D, $K'$
generates an indiscompositive property $K'_{\leq}$; in fact each graph
in $K'_\leq$ consists of a clique, isolated vertices,
and some end-vertices adjacent to some distinct vertices of the clique.  

Every indiscompositive or additive induced-hereditary property is induced-hereditary compositive and, therefore, an intersection property.  Is there a property of the sets $\cF$ that determines whether the intersection property is indiscompositive?  additive?

\section{Factorisation of additive hereditary properties\label{sec-alt-canonical}}

A property is {\em hereditary\/} if it is closed under subgraphs.  A property $\cP$ is {\em hereditary compositive\/} if it is hereditary and, for any $G$ and $H$ in $\cP$, there is a graph in $\cP$ containing both $G$ and $H$ as subgraphs.  Hereditary compositive properties can be characterised in essentially the same way as 
in parts A through D of Theorem~\ref{sch}, with appropriate changes in the definitions. 
We use this characterisation here to give a new proof of the result of Mih\'ok et al.~\cite{uft-1} that hereditary compositive and additive hereditary properties have a factorisation into properties with the same restrictions.  We are grateful to Jim Geelen for suggesting this approach.

A set $\cG$ is a {\em $\subseteq$-generating set\/} for $\cP$ if $\lb \cG \rb = \cP$, i.e. $\cG \subseteq \cP$ and every graph in $\cP$ is a subgraph of some $G \in \cG$; it is ordered if its elements can be listed as $G_1 \subseteq G_2 \subseteq \cdots.$ 
The set $\cG[L]_{\subseteq}$ is  $\{G \in \cG \mid L \subseteq G \}$, and $H$ is a {\em $\subseteq$-generating graph\/} for $\cP$ if $\cP = \{G \mid G \subseteq H,\, |V(G)| \textrm{ finite}\}$. If $G$ is a \emph{proper} in\-du\-ced-sub\-graph of $H$, we write $G<H$. A graph $G$ is 
{\em $\cP$-maximal\/} if $G \in \cP$, but for all $e \notin G$, $G+e \notin \cP$.
The set of $\cP$-maximal graphs is $\cM(\cP)$.

\statementnoboxnumbered{Theorem\label{sch-her}}
{For any property $\cP$ the following are equivalent:
\begin{itemize}
\item[A.] $\cP$ is hereditary compositive;
\item[B.] $\cP$ has a $\subseteq$-generating set; moreover, for any graph $L \in \cP$, and any generating set $\cG$ of $\cP$, $\cG[L]_{\subseteq}$ is also a $\subseteq$-generating set;
\item[C.] $\cP$ has an ordered $\subseteq$-generating set $G_1 \subseteq G_2 \subseteq \cdots;$
\item[D.] $\cP$ has an ordered $\subseteq$-generating set $H_1 < H_2 < \cdots$ where each $H_i$ is $\cP$-maximal; and
\item[E.] $\cP$ has a $\subseteq$-generating graph $H$.\hfill$\square$
\end{itemize}
}
We leave the proof of Theorem \ref{sch-her} (especially C $\Rightarrow$ D) to the reader, and note that the A $\Leftrightarrow$ C is contained in \cite[Thm. 2.3]{bebrmo}.  We also leave to the reader the statement and proof of the additive hereditary
version. 

To get the factorisation of hereditary compositive and additive hereditary properties, we need the concept of decomposability.  (Since it is rather more complicated, we do not discuss here the analogous concept for indiscompositive properties.  However, it is in \cite{uft-2,discussiones}, and the indiscompositive factorisation is given in the manner of this paper in \cite{farr-thesis}.)
A graph $G$ is {\em decomposable\/} if it is the join of
two graphs; otherwise, $G$ is {\em indecomposable}. It is easy to see that $G$ is decomposable if and only if its complement $\overline{G}$ is disconnected; $G$ is the join of the complements of the components of $\overline{G}$, so every decomposable graph can be expressed uniquely as the join of indecomposable subgraphs, the 
{\em ind-parts\/} of $G$. The number of ind-parts of $G$ is the {\em decomposability number\/} $dc(G)$ of $G$.

If $\cP$ is hereditary and contains $K_r$, then it contains all graphs on at most $r$ vertices. Since $\cP$ does not contain all graphs, there is an integer $c = c(\cP)$ such that $K_c \in \cP$, but $K_{c+1} \notin \cP$. The set of $\cP$-maximal graphs with at least $c(\cP)$ vertices is $\cM^*(\cP)$; these are precisely the $\cP$-maximal graphs $G$ such that $K_1 + G \not\in \cP$.  We require the following two elementary results from \cite{uft-1}.

\stmtnum[\ \cite{uft-1}]{Lemma\label{max-char}}
{Let $\cP = \cP_1 \circ \cdots \circ \cP_m$, where the $\cP_i$'s 
are hereditary graph properties.
A graph $G$ belongs to $\cM(\cP)$ if and only if, for every $(\cP_1, \ldots,
\cP_m)$-partition $(V_1, \ldots, V_m)$ of $G$, the following holds: 
$G[V_i] \in \cM(\cP_i)$ for $i = 1,\ldots, m$, and $G = G[V_1] +
\cdots + G[V_m]$. 
Moreover, if $G \in \cM^*(\cP)$, then each $G[V_i]$ is in $\cM^*(\cP_i)$ and, thus, non-null; therefore, $m \leq dc(G)$.
}
\stmtnumend

\stmtnum[\ \cite{uft-1}\label{uft-1-9}]{Lemma\label{gen-0}}
{Let $\cP$ be a hereditary property and let 
$G \in \cM^*(\cP)$, $H \in \cP$. If $G \subseteq H$, then $dc(H) \leq dc(G)$. If we have equality, with $G = G_1 + \cdots + G_n$ and $H = H_1 + \cdots + H_n$ being the respective expressions as joins of ind-parts, then 
we can relabel the ind-parts of $H$ so that, for each $i$, $G_i \leq H_i$. 
}
\stmtnumend
\newline

The {\em decomposability number\/}
of a set $\cG$ of graphs is $dc(\cG):= \min\{dc(G) \mid G \in \cG\}$. Clearly, for any hereditary property $\cP$, $\cM^*(\cP)$ generates $\cP$; moreover~\cite{uft-1}, if $\cG$ generates $\cP$, then $dc(\cG) \leq dc(\cM^*(\cP))$, with equality if $\cG \subseteq \cM^*(\cP)$. With this motivation, the {\em decomposability number\/} of $\cP$ is defined as $dc(\cP) := dc(\cM^*(\cP))$. If $dc(\cP)=1$, then $\cP$ is {\em indecomposable\/}.  Note that $dc(\cP\circ \cQ)\ge dc(\cP)+dc(\cQ)$, so any factorisation of $\cP$ into hereditary factors has at most $dc(\cP)$ factors.
We can now give our new proof of the factorisation theorem of Mih\'ok et al, also extending the result to hereditary compositive properties.

\resnum[\ \cite{uft-1}\label{uft-1-15}]{Theorem\label{factorn-new}}
{A hereditary compositive property $\cP$ has a factorisation into $dc(\cP)$ 
(necessarily indecomposable) hereditary compositive factors. Moreover, when $\cP$ is 
additive, the factors can be taken to be additive too.
}
Since $\cP$ is compositive, by Theorem~\ref{sch-her} it has an ordered generating set 
of $\cP$-maximal graphs, say $H_1 < H_2 < \cdots.$
Fix some graph $J \in \cM^*(\cP)$ with decomposability $dc(\cP)$; by omitting the 
first few $H_i$'s, we can assume that $J$ is contained in $H_1$. Then, by 
Lemma~\ref{gen-0}, each $H_i$ has exactly $d := dc(\cP)$ ind-parts, which we can 
label $H_{1,i}, \ldots, H_{d,i}$, so that, for $j = 1, \ldots, d$, $H_{j,1} \leq H_{j,2} 
\leq H_{j,3} \leq \cdots.$ Let $\cP_j$ be the hereditary property generated by the 
$H_{j,i}$'s; note that $\cP_j$ is compositive because the $H_{j,i}$'s are ordered. We 
claim that $\cP = \cP_1 \circ \cdots \circ \cP_d$.

If $G$ is in $\cP_1 \circ \cdots \circ \cP_d$, then $V(G)$ has a partition $\{V_1, 
\ldots, V_d\}$ such that, for each $j$, $G[V_j] \in \cP_j$. So there exist $k_1, \ldots, 
k_d$ such that, for each $j$, $G[V_j] \subseteq H_{j,k_j}$. Taking $k := \max\{k_1, 
\ldots, k_d\}$ we have $G \subseteq H_{1,k} + \cdots + H_{d,k} = H_k$, so $G$ is in $\cP$. Conversely, if $G$ is in $\cP$, it is contained in some $H_k$, and it is easy to 
see that it has a $(\cP_1, \ldots, \cP_d)$-partition.
\newline

If $\cP$ is additive, we claim that each $\cP_j$ is also additive. For each $1 \leq r,s 
\leq d$ such that $\cP_r \setminus \cP_s \not= \emptyset$, fix a graph $X_{r,s} \in 
(\cP_r \setminus \cP_s)$. By omitting some graphs from our ordered generating set, we 
can assume that $X_{r,s} \subseteq H_{r,1}$ for each $r$ and $s$. 

To prove additivity of $\cP_j$ it is sufficient to show that, for all $i$, $2H_{j,i}$ is 
contained in some $H_{j,i'}$. By additivity of $\cP$, $(d!+1)H_{i}$ is contained in 
some $H_{i'}$. By Lemma~\ref{gen-0}, for each copy of $H_{i}$, there is a 
permutation $\phi$ of $\{1, \ldots, d\}$ such that, for each $k$, $H_{k,i} \leq 
H_{\phi(k),i'}$. Since there are only $d!$ possible permutations, there are two copies 
of $H_{i}$ with the same permutation, so for some $\phi$ we actually have, for each 
$k$, $2H_{k,i} \leq H_{\phi(k),i'}$. 

Now, $\cP_j \subseteq \cP_{\phi(j)}$ (otherwise $X_{j,\phi(j)} \notin \cP_{\phi(j)}$ is 
contained in $H_{j,i} \subseteq H_{\phi(j),i'}$, a contradiction). If $\phi^t(j) = j$, 
then by repeating this argument, we get $\cP_j \subseteq \cP_{\phi(j)} \subseteq
\cP_{\phi^2(j)} \subseteq \cdots \subseteq \cP_{\phi^t(j)} = \cP_j$. So $2H_{j,i} \leq H_{\phi(j),i'} \in \cP_{\phi(j)} = \cP_j$.
%
\resnumend
\newline

Some comments are in order here. Dorfling \cite{mjd-thesis} proves (in a very elementary way) that if $\cP$ is an additive hereditary property that is not irreducible, then $\cP$ has a factorisation into two additive hereditary factors.  From this it is a triviality that every additive hereditary property has a factorisation into irreducible additive hereditary factors.  Theorem \ref{factorn-new} shows that every additive hereditary property $\cP$ has a factorisation into $dc(\cP)$ {\em indecomposable\/} additive hereditary factors.   This shows the important fact that an irreducible property is also indecomposable; the converse is trivial.  The proof of the uniqueness of the factorisation involves showing there is at most one factorisation into {\em indecomposable\/} factors.  Theorem \ref{factorn-new} is required then to deduce that there is a unique factorisation into {\em irreducible\/} factors.

\section{Infinite graphs\label{sec-infinite}}
In this section we discuss which parts of Theorems~\ref{sch}--\ref{sch-her} still hold when we allow infinite graphs. 
If the graphs have at most $\kappa$ vertices, for some infinite cardinal $\kappa$, then a generating graph for $\cP$ is now a graph $H$ such that $\cP =\{G \leq H \mid |V(G)| \leq \kappa\}$, while the intersection property generated by a family $\cF$ is $\cI(\cF) := \{G(\cF') \mid \cF' \subseteq \cF, |\cF'| \leq \kappa\}$.
Besides, an ordered generating set now need not be countable.
Note that in\-du\-ced-he\-re\-di\-ta\-ry properties are still precisely the ones that have generating sets.

In Theorem~\ref{sch} we have
E $\Leftrightarrow$ D $\Rightarrow$ B $\Leftrightarrow$ A $\Leftarrow$ C
(for E $\Rightarrow$ D, take $H = G(\cF)$).
To see that (B,C) $\nRightarrow$ D, consider the property 
\[ \cP_{\infty} = \{k_1 \overleftrightarrow{P} \cup k_2 \overrightarrow{P} \cup L \mid 0 \leq k_1, k_2 < \aleph_0, L \in \cL\},\]
where $\overrightarrow{P}$ and $\overleftrightarrow{P}$ are the one-way and two-way countably infinite path, respectively, and $\cL$ is the property containing disjoint unions of at most countably many finite paths. 
This property has an ordered generating set $\overleftrightarrow{P} \leq 2 \overleftrightarrow{P} \leq \cdots.$
However, a generating graph for $\cP_{\infty}$ must contain $\aleph_0 \overleftrightarrow{P}$, and this must then be a graph in $\cP_\infty$, a contradiction.

In Theorems~\ref{indiscomp-character} and~\ref{additive-character} we have (C,D) $\Rightarrow$ B $\Leftrightarrow$ A, but (B,C) $\nRightarrow$ D.
Similarly, in Theorem~\ref{sch-her} we have D $\Rightarrow$ C $\Rightarrow$ B $\Leftrightarrow$ A $\Leftarrow$ E, but (B,C,D) $\nRightarrow$ E. Moreover, (B,C) $\nRightarrow$ D: if each graph in $\cQ$ is the union of 
a finite star and a (finite or infinite) number of  isolated vertices, 
the only $\cQ$-maximal graphs are the finite stars, which do not generate $\cQ$.

All other implications are open.

\section{Minimal forbidden (induced-)subgraphs\label{sec-her-ind-her}}
We have seen that the various compositive properties are characterised ``from above'' by generating sets and generating graphs. Hereditary and induced-hereditary properties are also characterised ``from below'', by excluding a list of subgraphs or 
induced-subgraphs. Hereditary properties are, \emph{a fortiori}, induced-hereditary, and so can be characterised in two ways; in this section we announce some straightforward results relating the two characterisations. 

A graph $H$ is a {\em minimal forbidden subgraph} for $\cP$ if $H \notin \cP$, but 
all the proper subgraphs of $H$ are in $\cP$. The set of minimal forbidden subgraphs of $\cP$ is $\Fs(\cP)$. Greenwell, Hemminger and Klerlein~\cite{greenwell} showed that $\cP$ is hereditary iff the graphs in $\cP$ are precisely those that have no subgraph in $\Fs(\cP)$. The set $\Fi(\cP)$ of minimal forbidden induced-subgraphs is defined similarly, and characterises induced-hereditary properties.

Incidentally, this is not always true if $\cP$ contains infinite graphs, as $\Fi(\cP)$ and $\cM(\cP)$ may be empty even if $\cP$ is hereditary. This happens, say, when $\cP$ is the property of having finitely many edges. Yet, as pointed out to us  by Jan Kratochv\'{\i}l, this property is characterised by forbidding the infinite clique, the infinite star and the infinite matching.
\newline

It is well known that additive (induced-)hereditary properties are precisely the ones whose minimal forbidden (induced-)subgraphs are connected.
However we have no significant results to determine, from the forbidden subgraphs, whether the property is hereditary compositive.

A hereditary property $\cP$ is in\-du\-ced-he\-re\-di\-ta\-ry, so it is characterised by both $\Fs(\cP)$ and $\Fi(\cP)$. Note that $\Fs(\cP) \subseteq \Fi(\cP)$.
How do we obtain $\Fs(\cP)$ from $\Fi(\cP)$, or vice-versa? When are the two sets equal? If $\cP$ is in\-du\-ced-he\-re\-di\-ta\-ry, can we recognise from $\Fi(\cP)$ whether $\cP$ is actually hereditary? We provide answers below.
\newline

For a set $S$ of graphs, let $\min_{\subseteq}S$ (and $\min_{\leq}S$) be the set of graphs in $S$ that have no proper (induced-)subgraph in $S$. Now, defining $\Fi(\cP, n) := \{G \in \Fi(\cP) \mid |V(G)| = n\}$, we have:
\begin{eqnarray*}
\Fs(\cP) &=& \min_{\subseteq}(\Fi(\cP)) = {\displaystyle \bigcup_n} \min_{\subseteq}(\Fi(\cP, n)), \textrm{ and} \\
\Fi(\cP) &=& \min_{\leq}\{H+ e_1 + \cdots + e_r \mid H \in \Fs(\cP), \{e_1, \ldots, e_r\} \subseteq E(\overline{H})\}.
\end{eqnarray*}
A {\em $\leq$-antichain} is a set of graphs, none of which contains another as an induced-subgraph; $\subseteq$-antichains are defined similarly. For a hereditary property, $\Fs(\cP)$ must be a $\subseteq$-antichain, and therefore a $\leq$-antichain, but $\Fi(\cP)$ need only be a $\leq$-antichain. In fact, we have the following:

\stmtnum[{\ \cite[Prop.\ 2.1.7]{farr-thesis}}]{Proposition\label{prop-1}}
{Let $\cP$ be hereditary. Then $\Fs(\cP)$ is finite if and only if $\Fi(\cP)$ is finite. $\Fs(\cP) = \Fi(\cP)$ if and only if $\Fi(\cP)$ is a $\subseteq$-antichain, if and only if (for each $n \in \mathbb{N}$) $\Fi(\cP, n)$ is a $\subseteq$-antichain.
}
\stmtnumend

\stmtnum[{\ \cite[Prop.\ 2.1.8]{farr-thesis}}]{Proposition\label{prop-2}}
{Let $\cP$ be hereditary. Then $\Fs(\cP) = \Fi(\cP)$ if and only if, for all $H \in \Fs(\cP)$, and for all $e \notin E(H)$, there is some graph $G \in \Fs(\cP)$ such that $G \leq H+e$.
}
%
\stmtnumend
\newline

There is a deceptively similar result that characterises those in\-du\-ced-he\-re\-di\-ta\-ry properties that happen to be hereditary.

\stmtnum[{\ \cite[Prop.\ 2.1.9]{farr-thesis}}]{Proposition\label{prop-3}}
{Let $\cP$ be an in\-du\-ced-he\-re\-di\-ta\-ry property. Then $\cP$ is hereditary if and only if, for all $H \in \Fi(\cP)$, and for all $e \notin E(H)$, there is some graph $G \in \Fi(\cP)$ such that $G \leq H+e$. In this case, $\cP = \{F \mid \forall G \in \Fi(\cP), G \nsubseteq F\}$.
}
%
%
\stmtnumend
\newline

As an example, consider the property $\cP$ of forests, where $\Fi(\cP)$ is the set of cycles. For any  graph $C_{k} \in \Fs(\cP)$, and any edge $e \notin C_k$, $C_k + e$ contains another cycle as an \emph{induced}-subgraph. Proposition~\ref{prop-3} then confirms that $\cP$ is hereditary; moreover, since no cycle contains another as a subgraph, Proposition~\ref{prop-1} guarantees that $\Fs(\cP)$ is also the set of cycles. Similar considerations hold when $\cP$ is the set of bipartite graphs.

When $\cP$ is the set of graphs with at most $k$ vertices, then $\Fi(\cP)$ contains all graphs on $k+1$ vertices. Clearly, if $G$ is in $\Fi(\cP)$, then $G+e$ is again in $\Fi(\cP)$, confirming that $\cP$ is hereditary. However, $\Fi(\cP)$ is definitely not a $\subseteq$-antichain, and in fact $\Fs(\cP) = \{\overline{K}_{k+1}\}$ is a proper subset of $\Fi(\cP)$. This also follows from the fact that $\overline{K}_{k+1} \nleq (\overline{K}_{k+1} + e)$.


\end{document}